%

Title: Iterated Class Forcing
Author: Sy Friedman, MIT
E-mail:sdf@math.mit.edu

\magnification=\magstep1
\documentstyle{amsppt}
\NoBlackBoxes
\catcode`\@=11

\def\leftrightarrowfill{$\m@th\mathord\leftarrow\mkern-6mu%
  \cleaders\hbox{$\mkern-2mu\mathord-\mkern-2mu$}\hfill
  \mkern-6mu\mathord\rightarrow$}

\atdef@?#1?#2?{\ampersand@\setbox\z@\hbox{$\ssize
 \;\;{#1}\;$}\setbox\@ne\hbox{$\ssize\;\;{#2}\;$}\setbox\tw@
 \hbox{$#2$}\ifCD@
 \global\bigaw@\minCDaw@\else\global\bigaw@\minaw@\fi
 \ifdim\wd\z@>\bigaw@\global\bigaw@\wd\z@\fi
 \ifdim\wd\@ne>\bigaw@\global\bigaw@\wd\@ne\fi
 \ifCD@\hskip.5em\fi
 \ifdim\wd\tw@>\z@
 \mathrel{\mathop{\hbox to\bigaw@{\leftrightarrowfill}}\limits^{#1}_{#2}}\else
 \mathrel{\mathop{\hbox to\bigaw@{\leftrightarrowfill}}\limits^{#1}}\fi
 \ifCD@\hskip.5em\fi\ampersand@}
\atdef@-#1-#2-{\ampersand@\setbox\z@\hbox{$\ssize  
 \;\;{#1}\;$}\setbox\@ne\hbox{$\ssize\;\;{#2}\;$}\setbox\tw@
 \hbox{$#2$}\ifCD@
 \global\bigaw@\minCDaw@\else\global\bigaw@\minaw@\fi
 \ifdim\wd\z@>\bigaw@\global\bigaw@\wd\z@\fi
 \ifdim\wd\@ne>\bigaw@\global\bigaw@\wd\@ne\fi
 \ifCD@\hskip.5em\fi
 \ifdim\wd\tw@>\z@
 \mathrel{\mathop{\raise.5ex\hbox to\bigaw@{\hrulefill}}\limits^{#1}_{#2}}
\else\mathrel{\mathop{\raise.5ex\hbox to\bigaw@{\hrulefill}}\limits^{#1}}\fi
 \ifCD@\hskip.5em\fi\ampersand@}

\def\hookrightarrowfill{$\m@th\mathord\lhook\mkern-3.1mu%
  \cleaders\hbox{$\mkern-2mu\mathord-\mkern-2mu$}\hfill
  \mkern-6mu\mathord\rightarrow$}
\atdef@(#1(#2({\ampersand@\setbox\z@\hbox{$\ssize
 \;\;{#1}\;$}\setbox\@ne\hbox{$\ssize\;\;{#2}\;$}\setbox\tw@
 \hbox{$#2$}\ifCD@
 \global\bigaw@\minCDaw@\else\global\bigaw@\minaw@\fi
 \ifdim\wd\z@>\bigaw@\global\bigaw@\wd\z@\fi
 \ifdim\wd\@ne>\bigaw@\global\bigaw@\wd\@ne\fi
 \ifCD@\hskip.5em\fi
 \ifdim\wd\tw@>\z@
 \mathrel{\mathop{\hbox to\bigaw@{\hookrightarrowfill}}\limits^{#1}_{#2}}\else
 \mathrel{\mathop{\hbox to\bigaw@{\hookrightarrowfill}}\limits^{#1}}\fi
 \ifCD@\hskip.5em\fi\ampersand@}

\parindent20\p@
\advance\captionwidth@1in

\font\chapheadfont@=cmr10 scaled\magstep2
\font\chapheadmathfont@=cmmi10 scaled\magstep2
\font\chapheadmsa=msam10 scaled\magstep2
\font\chapheadmsb=msbm10 scaled\magstep2
\font\chapheadscriptmathfont=cmmi7 scaled\magstep2
\define\chapheadfont{\chapheadfont@\textfont1=\chapheadmathfont@%
	\scriptfont1=\chapheadscriptmathfont%
	\textfont\msafam=\chapheadmsa\textfont\msbfam=\chapheadmsb}
\outer\def\chapheading{\newpage
  \begingroup\raggedcenter@\interlinepenalty\@M \let\\\linebreak
  \chapheadfont\noindent\ignorespaces}

\newif\ifdatver
\define\datver#1{\ifdatver\global\def\datverp{}\else\datvertrue%
  \global\def\datverp{\par\boxed{\boxed{\text{Version: #1; Run: \today}}}}\fi}

\newif\ifbigdoc

\let\section\relax

\def\nosection{\newcodes@\endlinechar=10 \sect@}
{\lccode`\!=`\\
\lowercase{\gdef\sect@#1^^J{\sect@@#1!section\sect@@@}%
\gdef\sect@@#1!section{\futurelet\next\sect@@@}%
\gdef\sect@@@#1\sect@@@{\ifx\next\sect@@@\let
\next=\sect@\else\def\next{\oldcodes@\endlinechar=`\^^M\relax}%
 \fi\next}}}
\catcode`\@=\active
\define\protag#1 #2{{\hbox{\rm\ignorespaces#1\unskip}}#2}
\define\theprotag#1 #2{#2 {\rm\ignorespaces#1\unskip}}
\define\exertag #1 #2{\demo{\hbox{\rm(\ignorespaces#1\unskip)} #2}}

\define\figtagg#1{Figure $#1$}

\NoRunningHeads
\CenteredTagsOnSplits
\NoBlackBoxes

\def\today{\ifcase\month\or
 January\or February\or March\or April\or May\or June\or
 July\or August\or September\or October\or November\or December\fi
 \ \space\number\day, \number\year}

\define\QQ{\Bbb Q}

\define\Dom{\operatorname{Dom}}

\define\OS {\negthinspace\smallsetminus\negthinspace 0}

\define\dCc#1{\dot{\Cal{C}}^{#1}_c}

\define\nov#1{{\frac{n}{#1}}}

\define\({\bigl(}
\define\){\bigr)}
\UseAMSsymbols

\hsize 6.0truein
\vsize 8.2truein
\hcorrection{.2truein}
\vcorrection{.2truein}
\loadbold

\define\card{\operatorname{card}}
\define\Card{\operatorname{Card}}

\define\cof{\operatorname{cof}}

\redefine\phi{\varphi}
\document
\baselineskip=15pt

\font\bigtenrm=cmr12 scaled\magstep2
\centerline
{\bigtenrm {Iterated Class Forcing}}

\vskip10pt

\font\bigtenrm=cmr10 scaled\magstep2
\centerline
{\bigtenrm{Sy D. Friedman}
\footnote"*"{Research supported by NSF Contract\# 9205530-DMA.}}

\centerline
{\bigtenrm {M.I.T.}}

\vskip20pt

\comment
lasteqno 1@0
\endcomment

In this paper we develop the notion of ``stratified'' class forcing and show
that this property both implies cofinality-preservation and is preserved by
iterations with the appropriate  support. Many Easton-style and Jensen-style
forcings are stratified, as are some more exotic forcings obtained by mixing
these types together (see Easton [70], section 36 of Jech [78], 
Beller-Jensen-Welch [82], Friedman [90]).

As a sample application, cofinalities are preserved by an iteration of length
ORD where at even stages $i$,  an Easton-style forcing adds a Cohen set to
regular cardinals $\ge$ card $(i)$, at odd stages $i+1$  the class added at
stage $i$ is coded by a subset of the least infinite regular cardinal $\ge$ 
card(i) via the techniques of Friedman [93] or Friedman [94], and for any 
regular $\kappa$,  $\{i|p(i)$ is nontrivial below $\kappa\}$  is a subset of 
$\kappa^+$ of size $<\kappa$, for each condition $p$  in the iteration.

Jensen coding as in Beller-Jensen-Welch [82] is not stratified but obeys a
related property, called $\Delta$-stratification, which is also preserved by
iterations with the appropriate (larger) support. As a sample application,
the original form of Jensen coding can be used in the iteration of the
preceding paragraph, provided the Cohen sets are added with full support at
successor cardinals only and the condition stated at the end of that paragraph
is imposed only at successor cardinals.

We now define stratification, in the language of G\"odel-Bernays class theory.

\vskip5pt

\flushpar
{\bf Definition} \ $P$ (partially ordered by $\le$) is {\it stratified} if
there is a class $A$  such that $V=L[A]$ has a $\langle V,A\rangle$-definable
well ordering and:

(a) \ $P$  and $\le$ are $A$-definable. A condition in $P$  is a function $p$
on an initial segment of Card $=\{0\}\cup$ Infinite Cardinals, where if $q$
extends $p$  as a function, $q(\gamma)=\emptyset$ for all
$\gamma\in\Dom(q)-\Dom(p)$, then we identify $p$  with $q$.  Also we require
that $p(\kappa)=\emptyset$ for singular $\kappa$ and the conditions with 
constant value $\emptyset$ 
are the weakest in $P\/$.  Lastly, $\{p|p\restriction\kappa^+\in
H_{\kappa^{+}}\}$ is dense for each $\kappa\in\Card$. 

(b) \ ($\kappa$-Density Reduction) \ Let $\kappa$ be regular and define
$p\le_\kappa q$ if $(p\le q\text{ and }
p\restriction\kappa^+=q\restriction\kappa^+)$.  Then $p\le
q\longrightarrow\exists r\le_\kappa  q\exists s\le
r(s\restriction[\kappa^+,\infty)=r\restriction[\kappa^+,\infty)$  and 
$s\le p)$  and $p\le q\longleftrightarrow\exists r\exists s(p\le r, p\le s$ 
and $r\restriction\kappa^+ = q\restriction\kappa^+$,
$s\restriction[\kappa^+,\infty) = q\restriction[\kappa^+,\infty))$.  
If $D$  is an $A$-definable dense class, $p$  a condition then $\exists q\le_\kappa
p\exists d\subseteq D(\card(d)\le\kappa$  and every $r\le q$ can be extended
to $s$  such that $s\restriction
[\kappa^+,\infty)=r\restriction[\kappa^+,\infty)$  and $s$ extends an element
of $d\/$).

(c) \ ($\kappa$-Definable Closure) \ For infinite regular $\kappa$  there are
$\underset\sim\to\Pi_n^A$ operators $F_n(x,\kappa,p)$ for $0<n\in\omega$ such
that $F_n(x,\kappa,p)\le_{\kappa} p$ for all $p$  and whenever 
$p_0\ge_\kappa p_1\ge_\kappa\dots$  is a $\Pi_n^A$ 
(in parameters from $\kappa\cup\{x\})$
sequence of length $\lambda\le\kappa$  such that for each $i<\lambda$,
$p_{i+1}\le_{\kappa_i}$ $F_n(x,\kappa_i,p'_i)$ for some $p'_i\le_{\kappa_i}p_i$
and regular $\kappa_i\ge\kappa$  then there is $p\le_\kappa p_i$ for all
$i<\lambda$.   

\proclaim{\bf Theorem 1} Suppose that $P$  is stratified. Then $P$  preserves ZFC
(relative to the class $A$  witnessing stratification), cofinalities and the
GCH.
\endproclaim

\demo{Proof} Using $\kappa$-Definable Closure and $\kappa$-Density Reduction
we get: If $\langle D_i|i<\kappa\rangle$  is a uniformly $A$-definable
sequence of dense classes and $p\in P$  then there is $q\le p$  and $d$
of cardinality $\le\kappa$  such that each $r\le q$ is compatible with an
element of $D_i\cap d$, for each $i$.  This implies that the forcing relation
for $\underset\sim\to{\Delta }_0^A$ sentences is $A$-definable and that
$A$-replacement, cofinalities are preserved.

To show that GCH is preserved let $p\Vdash\tau\subseteq\kappa$,  where $\kappa$ is
an infinite cardinal and define $D_i=\{q|q$ is incompatible with $p$  or
$q\Vdash i\in\tau$ or $q\Vdash i\notin\tau\}$.  Then $D_i$ is dense and
$A$-definable, uniformly in $i<\kappa$.  First suppose that $\kappa$ is
regular. Then by $\kappa$-Density Reduction and $\kappa$-Definable Closure
there is $q\le p$ and $d$ of size $\le\kappa$  such that $r\in
d\longrightarrow
r\restriction[\kappa^+,\infty)=q\restriction[\kappa^+,\infty)$,
$r\restriction\kappa^+\in H_{\kappa^{+}}$ and  $\{r\in d|r\Vdash i\in\tau$ or
$r\Vdash i\notin\tau\}$ is predense below $q$ for each $i<\kappa$.  ($X$ is
predense below $q$ if every extension of $q$ is compatible with an element of
$X.)$ Thus in the generic extension each subset of $\kappa$  is determined by
a $\kappa$-sequence of size $\le\kappa$ subsets of $H_{\kappa^{+}}$ of the
ground model, so the property $2^\kappa=\kappa^+$ is preserved.

When $\kappa$  is singular the above argument can be repeated, using
$\kappa_i$-Density Reduction and cof$(\kappa)$-Definable closure, where
$\langle\kappa_i|i<\cof(\kappa)\rangle$ is a sequence  of regular cardinals
converging to $\kappa$. \hfill{$\dashv$ }
\enddemo

To preserve stratification under iteration we must discuss strong witnesses
and diagonal supports. 

\flushpar
{\bf Definition} \ $P\/$ is stratified with {\it strong witness\/} $A\/$ if in the
definition of stratified, $V=L[A]$ has a $\underset\sim\to{\Delta}^A_1$-well 
ordering,  $P,\le,\Card$ are $\underset\sim\to\Delta_1^A$ and there
is a $\underset\sim\to\Sigma_1^A$ function $f(x,\kappa,p)=(q,d)$ such that if
$x$ is an index for a $\underset\sim\to\Sigma_1^A$ dense class $D,p$ a
condition and $\kappa$ infinite and regular then $q\le_\kappa p$, card $(d)\le
\kappa, d\subseteq D$ and every $r\le q$ can be extended to $s$ such that
$s\restriction [\kappa^+,\infty)=r\restriction[\kappa^+,\infty)$  and $s$
extends an element of $d$. 

\proclaim{\bf Proposition 2} If $P$ is stratified then $P$ has a strong witness.
\endproclaim

\demo{Proof} Suppose $A$  witnesses that $P$ is stratified. Then we can choose
$A^*=\underset\sim\to\Sigma^A_N$ satisfaction, for a large $N$  so that $V$
has a $\underset\sim\to\Delta _1^{A^*}$  well ordering and
$P,\le,\Card,\{(q,d,\kappa)\big|\card(d)\le\kappa$ and every $r\le q$ can be
extended to $s$  such that
$r\restriction[\kappa^+,\infty)=s\restriction[\kappa^+,\infty)$  and $s$
extends an element of $d\}$ are $\underset\sim\to\Delta _1^{A^*}$. Then the
desired $\underset\sim\to\Sigma_1^A$ function $f(x,\kappa,p)$ exists. Finally
define the $\underset\sim\to\Pi_n^{A^*}$ operator $F^*_n(x,\kappa,p)$ to be
$F_{N+n}(x,\kappa ,p)$ where $\langle
F_m(x,\kappa,p)\big|0<m\in\omega\rangle$ comes from $A$. \hfill{$\dashv$}
  
\enddemo

Strong witnesses help us control the definability of the forcing relation.

\proclaim{\bf Theorem 3} Suppose that $A$  is a strong witness to the
stratification of $P$.  Then the forcing relation for $P$  restricted to
$\underset\sim\to\Sigma_1^A$  sentences is {\it densely-}$\Sigma_1^A$:
there is a  $\underset\sim\to\Sigma_1^A$ relation $p\overset{*}\to\Vdash\phi$
$(p\in P,\phi 
\underset\sim\to\Sigma_1^A)$   such that
$p\overset{*}\to\Vdash\phi\longrightarrow p\Vdash\phi$  and
$p\Vdash\phi\longrightarrow\exists q\le p$, $q\overset{*}\to\Vdash\phi$. 
\endproclaim

\demo{Proof} It suffices to prove this for $\underset\sim\to\Delta _0\phi$, 
by looking at witnesses.  We show by a $\Sigma_1^A$ induction on $\phi$ that
given $p$  we can (in a $\Sigma_1^A$ way) find $q\le p$ and $i\in\{0,1\}$ such
that either $q\Vdash\phi$, $i=0$  or $q\Vdash\neg\phi,i=1$.  This will prove
the 
Theorem since we can then take $q\overset{*}\to\Vdash\phi\longleftrightarrow$
For some $p,(q,0)$ arises from $p,\phi$ as above. 

The interesting case of the induction is the bounded quantifier: \ Suppose
$\phi$  is $\ \forall\ x\in\sigma\psi(x)$ where $\sigma$ is a term of rank
$\alpha$.  By induction we can effectively extend $p$ to decide any instance
$\psi(\tau)$, rank$(\tau)<\alpha$. If one of these extensions $q_\tau$ forces
$\tau\in\sigma\wedge\neg\psi (\tau )$ then we can take  the desired $q$ to be
$q_\tau\Vdash\neg\phi$.  Otherwise, we can build uniformly
$\underset\sim\to\Sigma_1^A$  dense classes $D_\tau$, rank$(\tau)<\alpha$ of
conditions forcing $\tau\in\sigma\longrightarrow\psi(\tau)$.  As $A$ is a {\it
strong} witness we can effectively find $d$  and $q\le p$  such that each
$D_\tau\cap d$ is predense below $q$,  where $d$ is a set. But this $q$ forces
$\ \forall\ x\in\sigma\psi(x)=\phi$. \hfill{$\dashv$

\enddemo

We are ready to discuss stratified iterations.

\vskip5pt

\flushpar
{\bf Definition} \ $\langle P_i|i<\alpha\rangle$ (where $P_{i+1}=P_i*\QQ_i$,
$P_\lambda\subseteq$ Inverse Limit $\langle P_i|i<\lambda\rangle$ for limit
$\lambda)$  is a {\it stratified iteration} if for some class $A\subseteq
ORD,A$  strongly witnesses that $P_0$ is stratified, for each
$i+1<\alpha,\emptyset_i\Vdash_i \QQ_i$ is stratified with strong witness
$\langle A,G_i\rangle$ via some $F^i_n,f^i$ $(\emptyset_i=$ weakest condition
in $P_i,\Vdash_i=$ forcing for $P_i,G_i=$ generic for $P_i)$  and $\QQ_i,f^i$ 
are $\underset\sim\to\Delta_1^{A,G_i}$ uniformly in $i<\alpha$, $F^i_n$ is 
$\underset\sim\to\Pi_n^{A,G_i}$ uniformly in $i<\alpha$, for each $n>0$.  Such
an iteration has {\it short diagonal supports} if for $p\in P_j$ and infinite
regular $\kappa$,  $\{i|i<j$ and $p\restriction i\not\Vdash_i$ $\ \forall\
\gamma<\kappa^+,p(i)(\gamma)=\emptyset\}$  is a subset of $\kappa^+$ of size
$<\kappa$ (and this is the only restriction on supports).

\proclaim{\bf Stratification Theorem} Suppose $\langle P_i|i<\alpha\rangle$ is a
stratified iteration with short diagonal supports and GCH holds. Then
$P_\alpha$ is isomorphic to a stratified forcing (definably relative to a
class $A$  witnessing stratification). 
\endproclaim

\demo{Proof} First we note that in the definition of stratified, we may assume
one further condition about the operators $F_n(x,\kappa,p):$
$\kappa_1\le\kappa_2$ both regular, $p(\gamma)=\emptyset$ for all
$\gamma<\kappa_2\longrightarrow F_n(x,\kappa_1,p)=F_n(x,\kappa_2,p)$.  For, we
may achieve this property by redefining $F_n$  to be
$F^*_n(x,\kappa,p)=F_n(x,\kappa(p),p)$ where $\kappa(p)=$ least $\gamma$ such
that $p(\gamma)\ne\emptyset$, if $\kappa\le\kappa(p)$ and
$F^*_n(x,\kappa,p)=F_n(x,\kappa,p)$  otherwise.

We prove the Main Theorem by induction on $\alpha$, maintaining the coherence
property that the isomorphism of $P_\alpha $ with a stratified forcing 
$P^*_\alpha$ extend the (inductively produced) isomorphism of $P_\beta$ with the 
stratified forcing $P^*_\beta$ for $\beta<\alpha$,  viewing $P_\beta$  as a 
subforcing of $P_\alpha$  in the natural way (and $P^*_\beta$ as a subforcing of
$P^*_\alpha$ in a natural way that will be evident from the  construction).
Also if $A$ is our given witness to the stratification of the iteration then
$A$ will serve as a strong witness to the stratification of each $P^*_\alpha$.

The result in vacuous for $\alpha=0$  or $1$.  Suppose that $\alpha=\beta+1$
is a successor ordinal $>1$.  By induction $P_\beta$  is isomorphic to a
stratified forcing $P^*_\beta$  and let $\le^\beta_\kappa$,
$F^\beta_n(x,\kappa,p^\beta)$  come from the stratification of $P^*_\beta$.
Also $\emptyset_\beta\Vdash_\beta \QQ^*_\beta$ is stratified 
($\emptyset_\beta=$  weakest condition of $P^*_\beta$, $\Vdash_\beta=$ forcing
for $P^*_\beta,\QQ^*_\beta=$ the $P^*_\beta$-name for $\QQ^\beta)$  and
let $\le_\kappa,F_n(x,\kappa,q)$ result from this.

By Theorem 3, $\Vdash_\beta$ is densely $\underset\sim\to\Sigma_1^A$ when
restricted  to $\underset\sim\to\Sigma_1^A$ sentences. By replacing $A$  by
$\underset\sim\to\Sigma_1^A$-Satisfaction, we may assume in fact that
$\Vdash_\beta$ is $\underset\sim\to\Sigma_1^A$ for $\underset\sim\to\Sigma_1$
sentences. It follows that  $\Vdash_\beta$ is $\underset\sim\to\Sigma_n^A$  
when restricted to $\underset\sim\to\Sigma_n$  sentences. 

Now we define $P^*_\alpha$ to essentially consist of all functions $f$  on an
initial segment of Card such that for some $p^\beta$ in $P^*_\beta$  and some
$q,p^\beta\Vdash_\beta q\in\QQ^*_\beta$ and for all $\kappa\in\Dom(f)$,
$f(\kappa)=\langle p^\beta(\kappa),q(\kappa)\rangle$ where $q(\kappa)$ is the
canonical term denoting the result of applying the function denoted by $q$ to
$\kappa$.  However we must make two small modifications: insist that if
$p^\beta(\kappa)=\emptyset$ and $p^\beta\Vdash_\beta q(\kappa)=\emptyset$ 
or undefined then $f(\kappa)=\emptyset$ (instead of $\langle\emptyset$, a term 
for $\emptyset\rangle)$; also insist that Dom$(f)$ contains Dom$(p^\beta)$ and
rank$(q)<\cup\Dom(f)$, so that $\emptyset_\beta$ $\Vdash_\beta$ 
Dom$(q)\subseteq \Dom(f)$. Then clearly $P^*_\alpha$ is isomorphic to
$P_\alpha$ when $P^*_\alpha$ is ordered in the natural way (by ordering the
corresponding pairs $\langle p^\beta,q\rangle$ in $P^*_\beta *\QQ^*_\beta)$.
It is easy to verify condition (a) and the first part of (b) in the definition
of stratification. 

Next we demonstrate $\kappa$-Density Reduction for $P^*_\alpha$.  For
notational purposes we think of a condition in $P^*_\alpha$ as an element of
$P^*_\beta *\QQ^*_\beta$ (isomorphic to $P^*_\alpha)$.  Suppose $D\subseteq
P^*_\alpha$ is dense and $A$-definable and $(p^\beta,q)\in P^*_\alpha$.
Consider $D^{G^*_\beta }=\{q_0\in\QQ^*_\beta|(p^\beta_0,q_0)\in D$ for some
$p^\beta_0\in G^*_\beta\}$ where $G^*_\beta$ denotes the $P^*_\beta$-generic.
Then $D^{G^{*}_{\beta }}\subseteq\QQ^*_\beta$ is forced by $\emptyset_\beta\in
P^*_\beta$ to be dense. So by $\kappa$-Density Reduction for $\QQ^*_\beta$,
$\emptyset^*_\beta$ also forces that $\{q_0\in\QQ^*_\beta$ 
$\big| q_0\/$ reduces $D^{G^*_\beta}\/$ below $\kappa^+\/$, to size 
$\le \kappa \}\/$ is $\le_\kappa\/$-dense on $\QQ^*_\beta\/$
(i.e., that for any $q_0$ there is $q_1\le_\kappa q_0$ such that for some 
$d\subseteq D^{G^{*}_{\beta }}$ of size $\le\kappa$, every $r\le q_1$ can be 
extended to $s$ such that $s,r$ agree $\ge\kappa^+$ and $s$ extends an element 
of $d)$.  Thus $D_\beta=\{p_0^\beta|$ For some $q_0,d_0,p_0^\beta\Vdash_\beta
q\notin\QQ^*_\beta$ or $q_0\le_\kappa  q$ reduces $D^{G^{*}_{\beta }}$ below
$\kappa$,  to $d_0$, card $(d_0)\le\kappa\}$ is dense on $P^*_\beta$. Let
$p_0^\beta\le^\beta_\kappa$ $p^\beta$  reduce $D_\beta$ below $\kappa^+$, to
size $\le\kappa$,  by $\kappa$-Density Reduction for $P^*_{\beta}$.

Then we can form terms $q_0,d_0$  such that $p_0^\beta\Vdash_\beta
q_o\le_\kappa q$, $q_0$ reduces $D^{G^*_\beta }$ below $\kappa^+$, to $d_0$ of
size $\le_\kappa$.  For each $i<\kappa$ it is dense below $p_0^\beta$ to force
some $p_0^\beta(i)\in G^*_\beta$,  $(p^\beta_0(i),d_0(i))\in D$,  where
$d_0(i)=i^{\text{th}}$ element of $d_0$.  Finally, by $\kappa$-Density
Reduction and $\kappa$-Definable Closure for $P^*_\beta$,  we can assume that
$p_0^\beta$ reduces all of these dense sets below $\kappa^+$, to size
$\le\kappa$ and hence $(p^\beta_0,q_0)$ reduces $D$ below $\kappa^+$, to size
$\le\kappa$. 

To complete the successor case we need to define the operators
$F^{\beta+1}_n(x,\kappa,(p^\beta,q))$  and verify condition (c) in the
definition of stratified. We set $F^{\beta+1}_n(x,\kappa,(p^\beta,q))=$
``least'' $(p^\beta_0,q_0)$  such that $p_0^\beta\le^\beta_\kappa
F_n^\beta(x,\kappa,p^\beta)$  and 
$p^\beta_0\Vdash_\beta F_n(x,\kappa,q)=q_0$. 
Note that the property of $(p_0,q_0)$ stated here is 
$\underset\sim\to\Sigma_n^A$  (in the other parameters), so we take ``least''
in the sense of $\underset\sim\to\Sigma_n^A$-uniformization, so that
$F_n^{\beta+1}(x,\kappa,p^{\beta+1})$ is $\underset\sim\to\Sigma_n^A$.
Of course we must show that such a $(p_0^\beta,q_0)$ exists. Note that it is a
dense property of $p_0^\beta$ to force a value for $F_n(x,\kappa,q)$.  By
$\kappa$-Density Reduction there is $p_0^\beta$ reducing this dense property
to a set, with $p_0^\beta\le^\beta_\kappa$ $F^\beta_n(x,\kappa,p^\beta)$.
Thus we can form a term $q_0$ such that $p^\beta_0\Vdash_\beta
F_n(x,\kappa,q)=q_0$. 

The $\kappa$-Definable Closure of $P^*_{\beta+1}$ follows from the
$\kappa$-Definable Closure of $P^*_\beta$ (relative to $A)$  and the
$\kappa$-Definable Closure of $\QQ^*_\beta$ (relative to $\langle
A,G^*_\beta\rangle)$.  Also $\le^{\beta+1}_\kappa$ is
$\underset\sim\to\Delta_1^A$-definable, uniformly in $\kappa$,  using the
facts that $\le^\beta_\kappa$ is uniformly $\underset\sim\to\Delta _1^A$,
$\le_\kappa$ is uniformly $\underset\sim\to\Delta _1^{A,G_\beta }$ and
the fact that $\Vdash_\beta$ is $\underset\sim\to\Sigma_1^A$  when restricted
to $\underset\sim\to\Sigma_1$ sentences.

Now we turn to the case where $\alpha$  is a limit ordinal. We take
$P^*_\alpha$ to consist of all functions $f$ on an initial segment of Card
such that for some $\langle f_\beta|\beta<\alpha\rangle$  in the inverse limit
of $\langle P^*_\beta|\beta<\alpha\rangle$ with short diagonal supports,
$f(\kappa)=\langle f_\beta(\kappa)|\beta<\alpha\rangle$ for all
$\kappa$  in Dom$(f)$; we also require that Dom$(f)\supseteq$ Dom$(f_\beta)$ for
each $\beta<\alpha$  and modify $f(\kappa)$ to be $\emptyset$ if
$f_\beta(\kappa)=\emptyset$ for all $\beta<\alpha$.  The $f$'s are ordered by
ordering the corresponding $\langle f_\beta|\beta<\alpha\rangle$'s.

We must show that $\{f\in P^*_\alpha|f\restriction\kappa^+\in H_{\kappa^+}\}$
is dense for each $\kappa$.  We actually show a bit more, for the purpose of
carrying out an inductive argument: if $\gamma<\kappa$ belong to Card,
$\gamma$ regular then $\{f\in P^*_\alpha|f\restriction[\gamma,\kappa ]\in
H_{\kappa^{+}}\}$ is $\le_\gamma$-dense (any $f$ can be $\le_\gamma$-extended
into this set; for $\gamma=0$ take $\le_\gamma=\le.)$  Note that this stronger
version follows from the weaker one, given $\gamma$-Density Reduction, so we
may inductively assume that it holds for $P^*_\beta,\beta<\alpha$.  Now we
induct on $\kappa$:  using short diagonal supports, we may assume that
cof$(\alpha)<\kappa$ as otherwise our given $f$  has the property that for
some $\beta_0<\alpha$, $f_\beta$ is the $\emptyset$-function below $\kappa^+$
for all $\beta_0\le\beta<\alpha$ (where $f$  comes from $\langle
f_\beta|\beta<\alpha\rangle)$ and so we can apply induction at $\beta_0$.  By
induction on $\kappa$  we can first extend $f$  to guarantee that
$f\restriction[\gamma,\cof(\alpha))$  belongs to $H_{\cof(\alpha)^+}$.  So we
may assume that $\gamma\ge\cof(\alpha)$.  Now, choose a cofinal
cof$(\alpha)$-sequence $\alpha_0<\alpha_1<\dots$ below $\alpha$  and
successively  extend $f=f_0$ to $f_1,f_2,\dots$ in cof$(\alpha)$ steps so
that $(f_{i+1}\restriction\alpha_i)$ on $[\gamma,\kappa ]$ belongs to
$H_{\kappa^+}$  and $f_{i+1}\restriction\alpha_j\ge_\gamma$
$F_1^{\alpha_j}(x,\gamma,f_i\restriction\alpha_j)$ for all $j\le i$, where
$F_1^{\alpha_i}$ comes from Definable Closure for $P_{\alpha_i}$ and
$x=\langle f,\gamma,\kappa,\langle\alpha_i|i < \cof\alpha\rangle\rangle$.  (We
abuse notation slightly; $f_i\restriction\alpha_i$  actually should be the
function $g(\delta)=f_i(\delta)\restriction\alpha_i.)$ Note that a simple
construction using the $F_1^{\alpha_j}$'s shows that $f_{i+1}$ as above does
exist. So we get that $f_\lambda$  is a condition for limit $\lambda$  and
$f_{\cof(\alpha)}$ is as desired.

If cof$(\alpha)\le\kappa$ we define $F^\alpha_n(x,\kappa,p)$ to be the least
$q\le_\kappa p$ such that $q\restriction\alpha_i\le_\kappa
F^{\alpha_{i}}_n(x,\kappa, p\restriction\alpha_i)$  for each $\alpha_i$ in a
fixed cof$(\alpha)$-sequence cofinal in $\alpha$.  If
$\kappa<\cof(\alpha)<\alpha$ and cof$(\alpha)$ is not the successor of a
regular cardinal then we obtain $F_n^\alpha(x,\kappa,p)$ by first choosing
$q_0\le_{\cof(\alpha)}p$ so that
$q_0\restriction\alpha_i\le_{\cof(\alpha)}F^{\alpha_i}_n(x,\cof(\alpha),
p\restriction\alpha_i)$  for each $\alpha_i$ 
and then $q_1\le_\kappa q_0$ so that $q_1\restriction\cof(\alpha)\le_\kappa
F_n^{\cof(\alpha)}(x,\kappa,q_0\restriction\cof(\alpha))$. If
$\kappa\le\lambda<\lambda^+=\cof(\alpha)<\alpha$ with $\lambda$  regular then 
we choose $q_0,q_1$ as above and then $q_2\le_\lambda q_1$ so that
$q_2\restriction\beta\le_\lambda F_n^\beta(x,\lambda,q_1\restriction\beta)$
for each $\beta$  such that $\cof(\alpha)\le\beta$  and
$q_1\restriction\beta\not\Vdash_\beta  q_1(\beta)$ is the $\emptyset$-function
below cof$(\alpha)$.  Finally if $\kappa<\alpha$ and $\alpha$ is regular then
choose $q_0$ as above $(\alpha_i=i)$  and then $q_1\le_\kappa q_0$ so that
$q_1\restriction\beta\le_\kappa F_n^\beta(x,\kappa,q_0\restriction\beta)$
where $\beta<\alpha$ is least so that $\beta\le\beta'\longrightarrow
q_0\restriction\beta'\Vdash_{\beta'}q_0(\beta')$ is the $\emptyset$-function
below $\alpha^+$.  Our construction guarantees that if
$q=F^\alpha_n(x,\kappa,p)$  and $\beta<\alpha$ then for some
$\kappa'\ge\kappa$ and $q',q\restriction\beta\le_\kappa
q'\le_{\kappa'}F^\beta_n(x,\kappa',p')$  for some $p'\le_\kappa
p\restriction\beta$.  The latter is used to verify $\kappa$-Definable Closure
when $\alpha$ is regular. (The other cases are straightforward, using our
extra hypothesis about the $F_n$'s stated at the start of the proof.)

Finally we must establish the second part of $\kappa$-Density Reduction for
$P^*_\alpha$.  (The first part is easy if $\cof(\alpha)\le\kappa$ and
otherwise follows inductively.) First suppose that $\alpha<\kappa^+$ and
choose an increasing cofinal $\alpha_0<\alpha_1<\dots$ of ordertype
cof$(\alpha)$.  Given $p\in P^*_\alpha$ and $A$-definable open dense $D$,  use
the $F^{\alpha_{i}}_n$  functions to successively extend
$p\restriction\alpha_i$ producing $q\le_\kappa p$ such that for each
$i<\cof\alpha$, $q\restriction\alpha_i$  reduces $D^{\alpha_{i}}$ below
$\kappa^+$ to size $\le\kappa$, where
$D^{\alpha_i}=\{r\restriction\alpha_i|r\in D\}$.  Now successively
$\le_\kappa$-extend $q=q_0$ to $q_1,q_2,\dots$ so that for each $\gamma$ there
is $x_{\gamma+1}$  defined on Card $\cap\kappa^+$ so that $x_{\gamma+1}\cup
q_{\gamma+1}$  is an element of $D$  yet is incompatible with each
$x_{\gamma'+1}\cup q_{\gamma'+1}$ for $\gamma'<\gamma$.  But for each $i$
there must be a stage $\gamma_i<\kappa^+$ such that for $\gamma\ge\gamma_i$,
$(x_{\gamma+1}\cup q_{\gamma+1})\restriction\alpha_i$ is compatible with some
$(x_{\gamma'+1}\cup q_{\gamma'+1})\restriction\alpha_i$ where
$\gamma'<\gamma_i$,  since $q\restriction\alpha_i$ reduces $D^{\alpha_i}$
below $\kappa^+$ to size $\le\kappa$.  Let
$\gamma=\cup\{\gamma_i|i<\cof\alpha\}<\kappa^+$.  Then $q_{\gamma+1}$ is
undefined so some $q_{\gamma'},\gamma'<\kappa^+$ reduces $D$  below
$\kappa^+$, to size $\le\kappa$. 

Now suppose that $\alpha\ge\kappa^+$.  We may assume that $\alpha=\kappa^+$ as
short diagonal supports requires that $p\in P^*_\alpha$  are trivial below
$\kappa^+$  on all but fewer than $\kappa$  coordinates, all below $\kappa^+$. 
But note that we can assume that conditions in $D$  when restricted to Card
$\cap\kappa^+$ belong to $H_{\kappa^+}$ and therefore can choose $q\le_\kappa
p$  and $\alpha_o<\kappa^+$ of cofinality $\kappa$  such that the conditions in
$D$  which are trivial below $\kappa^+$ on coordinates $\ge\alpha_0$ form a
set predense below $q$.  If we extend $q$ to $q_0\le_\kappa q$  reducing
$D^{\alpha_0}$  below $\kappa^+$ to size $\le\kappa$,  then $q_0$  in fact
reduces $D$  below $\kappa^+$ to size $\le\kappa$.  \hfill{$\dashv$ }
\enddemo

There are some important examples of cofinality-preserving class forcings that
are not stratified. Instead they may obey $\Delta$-stratification, which we
now consider.

\vskip5pt

\flushpar
{\bf Definition} \ $P$  is $\Delta$-{\it stratified} if it obeys the
definition of stratified where (b), (c) are restricted to successor cardinals
and in addition: whenever $0<n\in\omega$ and $\kappa$ is inaccessible, 
$p_0\ge p_1\ge\dots$ a $\Pi_n^A$ (in parameters from $\lambda\cup\{x\})$ sequence of 
length $\lambda\le\kappa$ and for each $i<\lambda$,
$p_{i+1}\le_{\kappa_{i}}F_n(x,\kappa_i,p'_i)$  for some
$p'_i\le_{\kappa_i}p_i$  and regular $\kappa_i\ge\aleph_{i+1}$, there is
$p\le_{\aleph_{i+1}}p_i$ for each $i$.  

$A$ is a {\it strong witness} to the $\Delta$-stratification of $P$  if it
obeys the definition of strong witness to stratification when restricted to
successor cardinals.  A $\Delta$-{\it stratified iteration} is just like a
stratified iteration but with stratified replaced by $\Delta$-stratified
everywhere.  Such an iteration $\langle P_i|i<\alpha\rangle$ has {\it
long diagonal supports} if for $p\in P_j$  and successor cardinals $\kappa$,
$\{i<j|p\restriction i\not\Vdash\ \forall\
\gamma\le\kappa,p(i)(\gamma)=\emptyset\}$ is a subset of $\kappa^+$ of size
$<\kappa$  and for inaccessible $\kappa\le j$,  $\{\bar\kappa<\kappa|$ For
some $\bar\kappa\le j'< j$, $p\restriction j'\not\Vdash_{j'}p(j')$ is
$\emptyset$ at $\bar\kappa\}$ is nonstationary in $\kappa$  (and these are the
only support restrictions). 

\proclaim{\bf Theorem 4} Suppose that $P$ is $\Delta$-stratified. Then $P$
preserves $ZFC$ (relative to the class $A$  witnessing
$\Delta$-stratification), cofinalities and the GCH. 
\endproclaim

\demo{Proof} As in the proof of Theorem 1, using $\Delta$-stratification at
$\kappa$  and stratification at $\bar\kappa^+<\kappa$, when $\kappa$ is
inaccessible. \hfill{$\dashv$ }

\enddemo

\flushpar
\proclaim{\bf $\Delta$-Stratification Theorem} Suppose $\langle
P_i|i<\alpha\rangle$ is a $\Delta$-stratified iteration with long diagonal
supports and GCH holds. Then $P_\alpha$  is isomorphic to a
$\Delta$-stratified forcing (definably relative to a class $A$  witnessing
$\Delta$-stratification). 
\endproclaim

\demo{Proof} We follow the proof of the Stratification Theorem. Note that
Theorem 3 still applies since its proof only uses that conditions (b), (c)
hold at cofinally many regular cardinals. We proceed by induction on $\alpha$.
For successor $\alpha$  our earlier proof still shows that (b), (c) hold at
successor cardinals. For $\Delta$-stratification at an inaccessible $\kappa$,
use $\Delta$-stratification for
$P^*_\beta(\alpha=\beta+1),\emptyset_\beta\Vdash_\beta$
$\Delta$-stratification for $\QQ^*_\beta$ to obtain $\bar
p\restriction\beta\le_{\aleph_{i+1}}p_i\restriction\beta$  for each $i$, $\bar
p\restriction\beta\Vdash_\beta$ there is $\bar
p(\beta)\le_{\aleph_{i+1}}p_i(\beta)$  for each $i$  and then $\le_{\kappa^+}$
extend $\bar p\restriction\beta$ to $p\restriction\beta$ so that for some term
$p(\beta)$, $p\restriction\beta\Vdash_\beta
p(\beta)\le_{\aleph_{i+1}}p_i(\beta)$ for each $i$,  using $\kappa^+$-Density
Reduction. Then $p\le_{\aleph_{i+1}}p_i$  for each $i$ is as desired. 

When $\alpha$ is a limit ordinal we define $P^*_\alpha$ as before and first
show that $\{f\in P^*_\alpha|f\restriction[\gamma,\kappa ]\in H_{\kappa^+}\}$
is $\le_\gamma$-dense for each successor $\gamma<\kappa,\gamma$ and $\kappa$
in Card. We do this by induction on $\kappa$,  noting that we may assume it
holds for $P^*_\beta,\beta<\alpha$.  Using (long) diagonal supports we may
assume that either $\alpha=\kappa$  is inaccessible or cof$(\alpha)<\kappa$.
If cof$(\alpha)$ is a successor or cof$(\alpha)^+<\kappa$  then the old
argument can be applied, using cof$(\alpha)$-Definable Closure or
cof$(\alpha)^+$-Definable Closure applied to $P^*_{\alpha_i}$,
$\alpha_i<\alpha$.  So we may assume that either $\alpha=\kappa$ is
inaccessible or cof$(\alpha)^+=\kappa$ where cof$(\alpha)$ is inaccessible. In
the latter case we choose a cofinal cof$(\alpha)$ sequence
$\alpha_0<\alpha_1\dots$  and successively $\le_\gamma$-extend our given
$f=f_0$ to $f_1,f_2,\dots$ in cof$(\alpha)$ steps so that
$(f_{i+1}\restriction\alpha_i)$ $(\kappa)\in H_{\kappa^{+}}$  and
$f_{i+1}\restriction\alpha_j\ge_{\aleph_{j+1}}$
$F_1^{\alpha_j}(x,\aleph_{j+1}\cup\gamma, f_i\restriction\alpha_j)$ for all
$j\le i$,  $x=\langle
f,\gamma,\kappa,\langle\alpha_i|i<\cof(\alpha)\rangle\rangle$.  Note
that by induction we may extend $g=f_{\cof(\alpha)}$ so that $g\restriction
[\gamma,\cof(\alpha)]\in H_\kappa$, as desired. Finally if $\alpha=\kappa$  is
inaccessible use Definable Closure to successively $\le_\gamma$-extend $f=f_0$
to $f_1,f_2,\dots$  in $\kappa$  steps choosing a continuous cofinal
$\kappa_0<\kappa_1<\dots$ below $\kappa$  such that
$f_{i+1}\restriction(\kappa_i,\kappa_{i+1})$ belongs to $H_{\kappa^+_{i+1}}$
and
$f_{i+1}(\kappa_i)=\emptyset$ for all $i$,  using the fact that
$\{\gamma<\kappa|f(\gamma)\ne\emptyset\}$ is nonstationary in $\kappa$.  Then
$f_\kappa$ is as desired.

If cof$(\alpha)\le\kappa$  or $\alpha$ is a successor cardinal nor
$\cof(\alpha)$ is neither inaccessible nor the successor of an inaccessible
then we define $F^\alpha_n(x,\kappa,p)$ as in the stratified case. If
cof$(\alpha)>\kappa$ is inaccessible then let $\alpha_0<\alpha_1<\dots$ be a
cofinal $\cof(\alpha)$-sequence  so that $\alpha_j\ge\aleph_{j+2}\cup\kappa$
for each $j<\cof(\alpha)$ and let $F^\alpha_n(x,\kappa,p)$ be a lower bound of
$p=p_0,p_1,\dots$ where $p_{j+1}$ is least so that
$p_{j+1}\restriction\alpha_{j'}$ 
$\le_{\kappa\cup\aleph_{j+1}}$
$F_n^{\alpha_{j'}}(x,\kappa\cup\aleph_{j+1},p_j\restriction
\alpha_{j'})$ for all $j'\le j$.  If
$\kappa\le\lambda<\lambda^+=\cof(\alpha)<\alpha,\lambda$ inaccessible then
similarly modify the earlier definition of $q_2$,  enumerating the relevant
$\beta$'s in $\lambda$  steps. 

$\kappa$-Density Reduction for successor $\kappa$  follows just as in the
stratified case. $\Delta$-stratification also follows as our construction
implies that if $p_{i+1}\le_{\aleph_{i+1}}$ $F^\alpha_n(x, \aleph_{i+1}, p_i)$ 
for  $i<\kappa$ $(\kappa$ inaccessible) then for cofinally many
$\alpha'<\alpha,p_{i+1}\restriction\alpha'\le_{\aleph_{i'+1}}$
$F_n^{\alpha'}(x,\aleph_{i'+1},p_i\restriction\alpha)$ for each $i$ (and
some $i'\ge i$ depending on $\alpha',i)$.   Also if $\beta<\alpha$ and
$p\le_{\kappa} q$ in $P^*_{\beta+1},p$ at $\beta=q$ at $\beta$  then
$F_n^{\beta+1}(x,\kappa,p)$ at $\beta$  equals $F_n^{\beta+1}(x,\kappa,q)$ at
$\beta$.  So given $p_0,p_1,\dots$ of length $\lambda$  as in the hypothesis
of $\Delta$-stratification at $\kappa$  for $P^*_\alpha$,  we can obtain the
desired lower bound $p$ by choosing $q\restriction\beta+1$  to be a lower
bound for $\langle p\restriction\beta\cup\{\langle\beta,p_i\text{ at
}\beta\rangle\}|i<\lambda\rangle$  and taking $p(\beta)=q(\beta)$. 
\hfill{$\dashv$ }

\enddemo

\flushpar
{\bf Examples}

(a) \ Jensen coding (Beller-Jensen-Welch [82]) is equivalent to a
$\Delta$-stratified forcing. It is dense to have $p(0)\ne\emptyset$ and
restricted to such conditions (together with the $\emptyset$ conditions)
condition (a) is satisfied. (We must reindex though: $p(\kappa)=p(\kappa^+)$
in Jensen's sense.) The first part of (b) is clear at successor $\kappa$  and
the second part is one of Jensen's key lemmas. For (c) we take
$F_n(x,\kappa,p)$ to be the least $q\le_\kappa p$ such that for
$\lambda\le\gamma\in\Dom(p),\gamma\in\Sigma^A_{n-1}\text{Hull
}(\gamma\cup\{x\})$,  $(q)_\gamma$  meets all predense $D$  on $P_\gamma$  in
$\Sigma^A_{n-1}\text{ Hull }(\gamma\cup\{x,p\})$.  Jensen's lemmas show that
such a $q$ exists and that (c) is satisfied (one can assume that all the
$\kappa_i$'s are equal by looking at their lim inf). The extra
$\Delta$-stratification condition also follows from Jensen's work.

(b) The modification of Jensen coding in Friedman [93] is equivalent to a
forcing that is both stratified and $\Delta$-stratified. It is densely
embeddable in the forcing defined in the same way (after reindexing) but where
at limit cardinals $\kappa$,  we allow $p\restriction\kappa$ to code only an
initial segment of $p_\kappa$  (and belong to the coding structure for that
initial segment). This allows one to prove (c) at inaccessibles. The thinning
that was done there in the limit coding enables one to prove (b) at
inaccessibles.

(c) The modification of Jensen coding in  Friedman [94] is stratified. The
proof of (b) at inaccessibles uses the fact that conditions have Easton
domains.

(d) Easton forcing (see Easton [70]) where a Cohen set is added to each
regular cardinal via an Easton product is stratified. (Take
$F_n(x,\kappa,p)=p$.) If, instead, the full product is used but only at
successor cardinals (no restriction on the domains of conditions) then
$\Delta$-stratification is obtained (but (b) will hold only at successors).
Without the restriction to successor cardinals one has a ``hybrid'' forcing
that is neither stratified nor $\Delta$-stratified. Iterating it would require
use of ``mixed support.''

(e) The forcing of Friedman [90] is a mixture of Jensen-style and Easton-style
forcing. It is equivalent to a stratified forcing, provided one of the
stratified modifications of Jensen coding (see (b), (c) above) is used.

(f) Backwards Easton forcings with Easton support (see Jech [78], section 36) 
are stratified provided at regular $\kappa$  one uses a $\kappa^+$-CC forcing 
of size $\le\kappa^+$.

\enddocument